\documentclass[11pt, letterpaper]{article}

\usepackage{hyperref}
\usepackage{authblk}
\usepackage{latexsym, amssymb, amsfonts,amsmath}
\usepackage{graphicx}
\usepackage[utf8]{inputenc}
\usepackage[english]{babel}
\usepackage[autostyle, english = american]{csquotes}
\usepackage[nottoc]{tocbibind}
\usepackage{upgreek}
\usepackage{stmaryrd}
\usepackage{qtree}
\usepackage{comment}
\usepackage{amsmath,amssymb,trimclip,adjustbox}
\usepackage{amsfonts}
\usepackage{mathtools}
\usepackage{etoolbox}
\usepackage{epigraph}
\makeatletter
    \newlength\epitextskip
    \pretocmd{\@epitext}{\em}{}{}
    \apptocmd{\@epitext}{\em}{}{}
    \patchcmd{\epigraph}{\@epitext{#1}\\}{\@epitext{#1}\\[\epitextskip]}{}{}
    \makeatother
\newcommand{\negr}[1]{\boldsymbol{#1}}

\newcommand{\LNTplus}{\mathcal{L}_\mathbb{NT}^+}
\newcommand{\Lplus}{\mathcal{L}^+}
\newcommand{\LNTtwo}{\mathcal{L}_\mathbb{NT}^{\textrm{(2)}}}
\newcommand{\ed}{\end{document}}
\newcommand{\PNP}{\textbf{P${\mathbf{=}}$?NP}}

\newtheorem{Theorem}{Theorem}

\newtheorem{Definition}[Theorem]{Definition}

\title{We must know --- We shall know}
\author{Jacob Zimbarg Sobrinho \\ email \href{mailto:jzimbarg@uol.com.br}{jzimbarg@uol.com.br} }
\affil{Instituto de Matemática e Estatística da USP}
\date{}

\begin{document}

\maketitle

\begin{center}
\sc To the memory of Carlos Edgard Harle\footnote{Carlos Edgard Harle was a Brazilian-German mathematician, deceased on January 16, 2020 at the age of 82. His dignified life has  been entirely devoted to academia: to research, student-advising and to the teaching of mathematics. He was the first mathematician to discover the concept of {\em isoparametric families of submanifolds}. On this occasion, it is an undisputable matter of  justice to honor him, both as scholar and man.}
\end{center}

\epigraph{Of his final endeavors in connection with mathematical knowledge, the last word has still not been spoken.
But when in this field a further development is possible, it will not bypass Hilbert but go through him.}
{\smallskip--- Arnold Sommerfeld, in his eulogy for David Hilbert\phantom{abcdef}}

\begin{abstract}
\noindent In this article, I focus on the resiliency of the {\PNP} problem. The main point to deal with is the change of the underlying logic from first to second-order logic. In this manner, after developing the initial steps of this change, I can hint that the solution goes in the direction of the coincidence of both classes, i.e., \textbf{P=NP}.
\end{abstract}

\section{Introduction}
In her wonderful biography of David Hilbert, Constance Reid \cite[p.~153]{Reid96}
mentions Hilbert's opinions about mathematical problems and mathematicians, by using the expression ``worthwhile problems." I would rather prefer the word ``resilient" instead of ``worthwhile,"  due perhaps to the subjectivity of the latter. The {\em Oxford English  Dictionary} defines the word ``resilient" as {\em  Resuming the original shape or position after being bent, compressed, or stretched.} In other words, ``resilient problems in mathematics" are those that resist attempts to their solution for a (relatively) long time, or have long proofs, which might turn them to be difficult to solve. In Gödel's article over the length of proofs \cite[p.~397]{Godel36}, he expressed the following statement:
\begin{quotation}
Thus, passing to the logic of the next higher order has the effect, not only of making provable certain propositions that were not  provable before, but also of making possible to shorten, by an extraordinary amount, infinitely many of the proofs already available.
\end{quotation}
Perhaps to make mathematical proofs shorter, it occurred to me trying to develop second-order logic, my intention being to discuss here the initial steps of such a development. It has been widely spread out in the foundational folklore that each standard model of set theory has its own second-order logic, and therefore, a universal system for it, as in  the first-order case, simply does not exist. Be as it may, the aim at developing second-order logic presupposes the philosophical commitment to an objective reality of sets (or some other entities) to be called {\em The Absolute}, and considered as the universe of {\em all} mathematical happenings. This has been eloquently conveyed as a ``genuine problem" by Professor Abraham Robinson
\cite[p.~556]{Robinson79}:

\begin{quotation}
The particular problem in the philosophy of mathematics that has the greatest fascination for me --- and I maintain that it is a genuine problem --- is that of the existence, or reality, or intelligibility, or objectivity, of infinite totalities.
\end{quotation}

The preceding considerations remind me the outstanding logician and Berkeley mathematics professor, Jack Howard Silver. It is widely known that Silver did have strong reservations as to the existence of measurable cardinals, despite the fact that his doctoral dissertation was fully based on Ramsey cardinals, an easy consequence of the existence of measurable cardinals. In a speech delivered by the great mathematician P.~J.~Cohen at Stanford University, he expressed his unbelief on the existence of inaccessible cardinals, for they  were too large.\footnote{For those interested in the subject of large cardinals, let me recommend Kanamori's masterpiece \cite{Kanamori09}.} Silver himself had reservations about the consistency of $\mathbb{ZF}$ set theory, and I caught somewhere in the internet remarks about his suspicions of third-order number theory. In that case, mathematics would be reduced simply to number theory and analysis. Since set theory still arises doubts in the minds of a few mathematicians, why not develop second-order logic in number theory? That's what I will do, and at the end I'll make some remarks of how these ideas can be extended to set theory.

\medskip
In September 1917, Hilbert visited Zürich where he delivered a talk which according to Reid \cite[p.~151]{Reid96} {\em $\ldots$ was devoted to a favorite subject --- the importance of the role of mathematics in the sciences --- and might have been entitled {\em ``In praise of the axiomatic method."}} There he brought up certain questions which I would rather call {\em Foundational Questions} that I mention in the sequel:

\medskip
\begin{itemize}
    \item[] {\sc FQ 1.} The problem of the solvability  in principle of every mathematical question.
    \item[] {\sc FQ 2.} The problem of finding a standard of simplicity for mathematical proofs.
    \item[] {\sc FQ 3.} The problem of the relation of content and formalism in mathematics.
    \item[] {\sc FQ 4.} The problem of the decidability of a mathematical question by a finite procedure.
\end{itemize}
\medskip
I dare saying that, as our science progresses, all of the foundational questions above will be fully and positively answered in the framework of second-order logic. To tackle FQ 1, I would like to rephrase it in the
following way: {\em The language of mathematics (in our case, number theory) is consistent and complete.} This brings us immediately to two stumbling blocks: Gödel's incompleteness theorems \cite{Godel31} and Lindström's model-theoretic characterization of first-order logic \cite{Vaananen12}.
Let me quote a small passage in Cohen's book \cite[p.~1]{Cohen66}:
`` $\ldots$ Gödel's incompleteness theorem
still represents the greatest obstacle to a satisfactory philosophy of
mathematics."

To circumvent Gödel's incompleteness theorems, I chose to extend further the language of first-order number theory, or Peano's arithmetic; let me substantiate this move by quoting Hilbert \cite[p.~177]{Reid96}:
\begin{quotation}
``Let us remember that {\em we are mathematicians} and that as mathematicians we have often been in precarious situations from which we have been rescued by the ingenious method of ideal ele\-ments $\ldots$. Similarly, to preserve the simple formal rules of aristotelian logic, we must {\em supplement the finitary statements with ideal statements.}"
\end{quotation}

\noindent And according to the concluding remark in Reid  \cite[p.~177]{Reid96}:
\begin{quotation}
Mathematics, under this view, would become a stock of two kinds of formulas: first, those to which the meaningful communications correspond and, secondly, other formulas which signify nothing but which are the ideal structure of the theory.
\end{quotation}

And she proceeds by quoting Hilbert again:
\begin{quotation}
``But in our general joy over this achievement, and in our particular joy over finding that indispensable tool, the logical calculus, already developed without any effort on our part, we must not forget the essential condition of the method of ideal elements --- {\em a proof of consistency.}"
\end{quotation}

As we see, the first-order language of Peano's arithemetic will be broadened to shelter the ideal elements of the new language. Another purpose is to prevent the definability of the new language within itself, thus precluding the derivation of the incompleteness theorems. To this end, I interpret the new formulas as having an infinitary (non-well-founded) tree-like format, whose content in the nodes is similar to the older first-order situation. If the newer language turns out to be {\em non}-definable within itself, we might then suppose that the corresponding axiomatic system is consistent and complete.
As to the non-validity of Lindström's theorem, my excuses will be lighter, by just mentioning that the emergent language is not a set, and therefore, external to set theory. Thus, Lindström's theorem does not apply. Each language of this class will be referred to as an {\em infinitary first-order number theoretic language.}\footnote{During the seventies, I was told that the great Alfred Tarski mentioned the following saying: {\em The infinitary sentences will constitute the future of logic.} But I doubt very much that Tarski had in mind the (infinitary) sentences that are being considered here.}

\section{Infinitary first-order languages}

Perhaps it would be easier to introduce infinitary first-order formulas by means of an example. Consider, for instance, the formula which expresses induction in number theory for $\phi(x)$:
        \begin{displaymath}
            (\phi(0)\land(\forall x(\phi(x)\rightarrow\phi(Sx)))\rightarrow(\forall x\phi(x))).
        \end{displaymath}

\noindent In polish notation this formula is somewhat more weird:
    \begin{displaymath}
        \rightarrow\land\phi0\forall x\rightarrow \phi x\phi Sx\forall x\phi x.
    \end{displaymath}
Finally, in tree form it acquires the following shape:

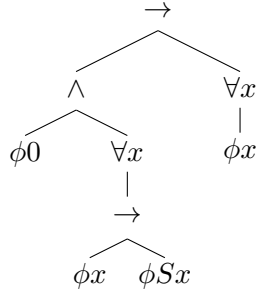
\Tree [.$\rightarrow$ [.$\land$ $\phi0$ [.$\forall x$ [.$\rightarrow$ $\phi x$ $\phi Sx$ ] ] ] [.$\forall x$ $\phi x$ ] ]

\medskip\noindent As we see, the notion of formula in tree-form can easily be extended to infinitary trees, which will turn out to be non-well-founded, i.e., they will admit infinite branches. The concept of free and bound occurrences of variables can be carried through straightforwardly to the new context as well as the substitution of variables, free or bound, without collisions. Also the notion of sub-formula is pretty much easy to infer. Now, it is time to give precise definitions:

\begin{Definition}
    An {\em infinitary first-order number theoretic formula} is given by a tree, (possibly) non-well-founded, satisfying the following conditions:
    \begin{enumerate}
        \item The total number of variables, free or bound, is finite;
        \item each node is either a boolean connective $\land$ (conjuction), $\lor$ (disjuction), $\lnot$ (negation), $\rightarrow$ (implication), $\leftrightarrow$ (equivalence); or a quantifier $\forall v_n$ (universal), $\exists v_n$ (existential);
        \item terminal nodes are atomic formulas of number theory $t_1=t_2$, $t_1<t_2$, where $t_1, t_2$ are terms built from variables $v_0, v_1,\cdots,v_n\cdots$, or from the constant $0$, through the application of the functional symbols $+$ (sum), $*$ (multiplication) or $S$ (successor); or may be given by symbols denoting continuing further trees.
    \end{enumerate}

    An {\em infinitary first-order number theoretic language} $\mathcal{L}^+$ is a {\em collection} of infinitary first-order number theoretic formulas satisfying the following conditions:
    \begin{enumerate}
        \item $\mathcal{L}^+$ includes all atomic formulas;
        \item $\mathcal{L}^+$ is closed under boolean combinations and first-order quantifications;
        \item $\mathcal{L}^+$ is closed under the operation of taking sub-formulas;
        \item $\mathcal{L}^+$ is closed under substitution of free or bound variables without collisions.
    \end{enumerate}
\end{Definition}

\noindent The reader should keep in mind the following observations:
\begin{itemize}
    \item the total number of variables in an infinitary first-order formula is finite, but not necessarily the number of their occurrences, which may be infinite;
    \item an infinitary first-order number theoretic language is not necessarily a set, but includes the first-order language for number theory, which is a set.
\end{itemize}

\medskip To finish this section, I will provide some notation to be used throughout the rest of the exposition: the first-order language for number theory will be denoted by $\mathcal{L}_{\mathbb{NT}}$; its second-order language, by $\mathcal{L}_{\mathbb{NT}}^{(2)}$; and a general infinitary first-order language for number theory will be designated by $\mathcal{L}^+$; I will characterize in a later section a special language for infinitary  first-order number theory, bridging the gap between first and second-order logic, that will be denoted by $\mathcal{L}_{\mathbb{NT}}^+$. Finally, I would point out that amongst the infinitary first-order number theoretic languages there are two which are extreme: a minimal, $\mathcal{L}_{\textrm{(inf)}}^+$, formed by short formulas only, and a maximal, $\mathcal{L}_{\textrm{(max)}}^+$, containing all infinitary formula-trees. Both are sets. A generic infinitary first-order number theoretic language $\mathcal{L}^+$ is in between the two, which entails that all formulas are sets, but not necessarily all languages $\mathcal{L}^+$:
\begin{displaymath}
    \mathcal{L}_{\textrm{(inf)}}^+\subseteq\mathcal{L}^+\subseteq\mathcal{L}_{\textrm{(max)}}^+.
\end{displaymath}

\section{Proofs and Decisions}
Before entering the subject of proofs and decisions, let me make a few remarks regarding the metamathematics of second-order logic. As said before, the metamathematical syntactical notions can be translated almost {\em verbatim} from first to second-order logic. But since in the latter we are dealing with infinitary formulas, a few adjustments are pertinent, mainly to maintain the finitary spirit of metamathematics.

\medskip First of all, it is always  possible to distinguish between object language and metalanguage proofs. The first are, so to speak, internal to the theory, and therefore, when we have in mind a model, the object language proofs are seated within. On the other hand, metalanguage proofs are always external, finitary and concrete, whereas in the first-order case, they might be infinitary and non-standard. In the framework of second-order logic, however, such a phenomenon should never happen, and proofs, metatheoretical or object language-like are always finitary or, to put it another way, isomorphic. Thus, for instance, were we able to succeed in showing within the second-order framework that $\mathbb{ZF}$ is first-order inconsistent, i.e., to derive the sentence $\lnot Con(\ulcorner\mathbb{ZF}\urcorner)$ in number theory, the internal proof of inconsistency could in principle be retrieved to a fully fledged proof of the metalanguage.

\smallskip
But, since we are dealing with infinitary objects in a proof, how is it that we are going to preserve finiteness? The answer to this question is simple: just consider the infinitary tree-formulas as potentially infinite objects. Since the proof-theoretical elements only depend on the initial stages of information in the trees, and not on the entire content of the formula-tree, this apparent limitation on the knowledge of formulas will not hurt the concept of proof. For instance, consider the logical axiom
$\Phi\rightarrow(\Psi\rightarrow\Phi)$. It can be easily seen that the knowledge of its being an axiom only depends on the 3 initial stages of its tree presentation, and not on the entire scope of its displayed information:

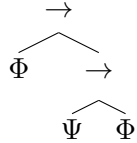
\Tree [.$\rightarrow$ $\Phi$ [.$\rightarrow$ $\Psi$ $\Phi$ ] ]

\medskip The notion of first-order proof is very well-known. For those interested in looking this concept up, I vividly recommend Barwise's essay \cite[p.~34--35]{Barwise77} on first-order logic. As said before, the concept of proof can be translated almost {\em verbatim} to infinitary first-order formulas. But I contend that in the second-order context, the notion of {\em decision} has to be introduced, and I wish to explain why. When in the context of first-order logic we prove a formula $\phi(v_0,\ldots v_k)$, what is really at stake is its universal closure, i.e., the sentence $\forall v_0\cdots\forall v_k\phi(v_0,\ldots,v_k)$. Whereas in second-order logic, we are interested in the validity of each instance $\phi(S^{n_0}0,\ldots,S^{n_k}0)$ obtained through assigning to each variable $v_i$ the term $S^{n_i}0$ in the whole decision, for $0\leq i\leq k$. In {\em Decision Theory} thus, the deductive rule of generalization cannot in general possibly apply.

I begin by quoting Matijasevich \cite[p.~38]{Matijasevich92}:
\begin{quotation}
The tenth problem is the only one of the 23 problems
that is (in today's terminology) a {\em decision problem}; i.e., a
problem consisting of infinitely many individual problems
each of which requires a definite answer: YES or
NO. The heart of a decision problem is the requirement
to find a single method that will give an answer
to any individual subproblem.
\end{quotation}

Nevertheless, despite the fact that Decision Theory has yet to be developed, I consider Proof theory and Decision Theory as sister theories. In fact, given an infinitary first-order number theoretic $\mathcal{L}^+$, a sentence $\Sigma\in\mathcal{L}^+$ admits a {\em proof} in $\mathcal{L}^+$ if there is a finite sequence of formulas $\pi^+\subset\mathcal{L}^+$ whose  last element is $\Sigma$ such that $\pi^+$ is a number theoretic proof (including or not some of the  Peano's axioms). On the other hand, a formula $\Phi(v_0,\ldots,v_k)$ is {\em decided} in $\mathcal{L}^+$ if there exist two finite sequences of formulas of $\mathcal{L}^+$, $\delta^+$ and $\delta^-$, whose last elements are $\Phi(v_0,\ldots,v_k)$ and $\lnot\Phi(v_0,\ldots,v_k)$ respectively, such that for every sequence $\langle S^{n_0}0,\ldots,S^{n_k}0\rangle$,

\[\delta^+(S^{n_0}0,\ldots,S^{n_k}0) \text{\ \ or\ \ } \delta^-(S^{n_0}0,\ldots,S^{n_k}0)\]

\noindent is a number theoretic proof. Here, $\delta^+(S^{n_0}0,\ldots,S^{n_k}0)$ (or its analogous counterpart $\delta^-(S^{n_0}0,\ldots,S^{n_k}0)$) is obtained by substituting all variables $v_i$ in $\delta^+$ (resp.~$\delta^-$) by $S^{n_i}0$. In case $\Sigma$ is {\em provable} in $\mathcal{L}^+$, we use the notation $\vdash_{\mathcal{L}^+}\Sigma$. In case $\Phi(v_0,\ldots,v_n)$ is {\em decided} in $\mathcal{L}^+$, we use the notation $\Vdash_{\mathcal{L}^+}\Phi(v_0,\ldots,v_n)$ (for decision theory reminds me subjectively of Cohen's notion of {\em forcing} in set theory). To finish this section, let me remark that Proof Theory is a subject linked to Foundational Question 2, whereas Decision Theory, to Foundational Question 4.

\section{The Absolute}
The Absolute for number theory will be constituted by two components: the set of natural numbers $\mathbb{N}$ and the predicates over numbers $\mathfrak{P}\mathbb{N}$. The Absolute $\mathfrak{A}$ will be, so to speak, the union of these two components:\footnote{If the reader wishes, she might refer to $\mathfrak{A}$ as a shortcut for {\em Analysis.}}
\begin{displaymath}
    \mathfrak{A}=\langle\mathbb{N},\mathfrak{P}\mathbb{N}\rangle = \mathbb{N}\cup\mathfrak{P}\mathbb{N}.
\end{displaymath}

\smallskip Predicates will be given by formulas of an infinitary language $\LNTplus$, to be characterized in the sequel, and are linked to the Foundational Question 3, whereas $\LNTplus$ is connected to the Foundational Question 1. In order to characterize $\LNTplus$, suppose we are given a numbering of all formulas of $\mathcal{L}_{\textrm{(max)}}^+=\langle\Phi_0,\Phi_1,\ldots,\Phi_n,\ldots\rangle.$  I'm going to consider a sequence $\langle\Gamma_n:n\varepsilon\mathbb{N}\rangle$ of infinitary first-order number theoretic languages satisfying the following conditions:
\begin{enumerate}
    \item $\Gamma_0$ contains all infinitary first-order number theoretic languages, including $\mathcal{L}_{\textrm{(max)}}^+$.
    \item If $\Lplus\in\Gamma_n$, in order for $\Lplus$ to remain in $\Gamma_{n+1}$, it is necessary that either $\Phi_n\notin\Lplus$ or $\Phi_n$ is decided in $\Lplus$, i.e., $\Vdash_{\mathcal{L}^+}\Phi_n$. Otherwise, $\Lplus$ is discarded when passing from $\Gamma_n$ to $\Gamma_{n+1}$.
    \item If $\Lplus\in\Gamma_n$ is inconsistent with number theory, i.e., there is a short sentence $\sigma=\Phi_n$ false in the standard model of number theory --- $\langle\mathbb{N},+,*,S,<,0\rangle\models\lnot\sigma$ --- then $\sigma$ is decidable favorably by $\Lplus$, i.e., $\Vdash_{\mathcal{L}^+}\sigma$,  and $\Lplus$ is discarded when passing from $\Gamma_n$ to $\Gamma_{n+1}$.
\end{enumerate}

\medskip It is fairly intuitive to check that any $\Lplus$ belonging to all of the $\Gamma_n$ is consistent, complete and all of its short sentences are true in the standard model of number theory. Besides, the sequence $\langle\Gamma_n:n\in\omega\rangle$ is decreasing, i.e.,
\begin{displaymath}
    \ldots\subseteq\Gamma_{n+1}\subseteq\Gamma_n\subseteq\ldots\subseteq\Gamma_1\subseteq\Gamma_0.
\end{displaymath}

\smallskip Now, I allow myself to state axiomatically the following principle:

\medskip
\noindent{\sc Monotonicity principle ---} There exists exactly one infinitary first-order number theoretic language $\LNTplus$ belonging to all of the $\Gamma_n$'s:
\begin{displaymath}
    \{\LNTplus\}=\bigcap_{n=0}^\infty\Gamma_n.
\end{displaymath}

Strictly speaking, I would have to consider the language $\LNTplus$ as given by a subsequence $\langle\Phi_{n_0},\ldots\Phi_{n_k},\ldots\rangle$. But as I am prone to commit abuses, I'll commit them once again and always write $\LNTplus = \langle\Phi_0,\ldots\Phi_n,\ldots\rangle$.

\smallskip The reader might wonder how is it possible that $\LNTplus$ is consistent and complete, since we are using only plain first-order logic in our deductions. The answer to this question lies in the fact that the used language is infinitary, and it might encode a lot of information about the universe in its structure. Let me give an example that might be considered as a plausibility argument. Consider the formula

\Tree [.$\forall x$ [.$\lor$ $x=0$ [.$\lor$ $x=S0$ [.$\lor$ $x=SS0$ [.$\lor$ $\vdots$ ] ] ] ] ]

\noindent I am not willing to say that this formula necessarily belongs to $\LNTplus$. But it is intuitive to check that it encloses a lot of information about $\mathbb{N}$.

\smallskip Now I'll come to grips with Foundational Question 3: my contention is that mathematical objects are self-referential predicates, i.e., predicates of a language whose predicates refer to themselves. In order to make this assertion clear, I initially must say what a predicate of a language is, and I restrict myself to predicates of $\LNTplus$. As we have seen, this language is supposed to be consistent and complete. Well, suppose we are given two formulas of $\LNTplus$, $\Phi(v_0,\ldots,v_n)$ and $\Psi(v_0,\ldots,v_n)$ (with the same number of free variables). Clearly, they are associated to two predicates. We say that $\Phi$ and $\Psi$ define the {\em same} predicte if they are equivalent with respect to the deductive system of $\LNTplus$:
\begin{displaymath}
    \vdash_\mathbb{NT}^+\forall v_0\ldots\forall v_n (\Phi(v_0,\cdots,v_n)\leftrightarrow\Psi(v_0,\ldots,v_n)).
\end{displaymath}
This definition was given by using proof theory, but it might have been given decision theoretically:
\begin{displaymath}
    \Vdash_\mathbb{NT}^+\Phi(v_0,\cdots,v_n)\Leftrightarrow\hskip3pt\Vdash_\mathbb{NT}^+\Psi(v_0,\cdots,v_n).
\end{displaymath}

In general, a number theoretic predicate is determined by a formula of $\LNTplus$. The predicate determined by the formula $\Phi(v_0,\ldots,v_n)$ will be denoted by
$$\llbracket\Phi(v_0,\ldots v_n) \rrbracket.$$
But sometimes I'll commit abuses by identifying a predicate to its defining formula:
$\llbracket\Phi(v_0,\ldots v_n) \rrbracket \approx\Phi(v_0,\ldots v_n).$ This relates our considerations to the Foundational Question 3.

\smallskip The definition of satisfaction and truth in the Absolute $\mathfrak{A}$ is a plain Tarski-style definition, individual variables ranging over numbers in $\mathbb{N}$, and predicate variables ranging over predicates over numbers in $\mathfrak{P}\mathbb{N}$. As I will need truth for existential second-order formulas only, I'll make it explicit: let $\exists X\varphi(X;v_0,\ldots v_n)$ an existential (short) formula of $\LNTtwo$. We will suppose, without loss of generality, that there is only one second-order quantification in the outer side of the formula. All of the others have been collapsed into a single one by means of suitable pairing functions. Then, in case $\mathfrak{A}\models\exists X\varphi(X;v_0,\ldots,v_n)[a_0/v_0,\ldots a_n/v_n]$ one can find a formula $\Phi(v_0,\ldots,v_n)\in\LNTplus$ satisfying the condition \[\vdash_\mathbb{NT}^+\varphi(X;v_0,\ldots,v_n)[\Phi/X;a_0/v_0,\ldots a_n/v_n].\]
$\Phi/X$ means that I have replaced  $X$ in $\varphi$ by the instances of $\Phi$ obtained through the substitution of $v_0,\ldots v_n$ by the corresponding variables in the occurrence of $X$ that is being replaced.

\section{\textbf{P=NP}} After having developed the {\em prolegomena} in previous sections, the second-order machinery is ready to be used to show the main consequence of this article, which is an answer to the \PNP\  problem. The argument will be short and easy. It uses one of the theorems of Descriptive Complexity, due to Ron Fagin
\cite[p.~3]{Halpern01}:

\begin{Theorem}\label{Fag} ({\sc Fagin}) A set of structures $\mathcal{T}$ is in \textbf{NP} iff there exists a second-order existential formula, $\exists X\varphi$
such that $\mathcal{T}=\{\mathcal{A}:\mathcal{A}\models\ulcorner\exists X\varphi\urcorner\}$. Formally, $\mathbf{NP} = SO\exists$.
\end{Theorem}

Before going into the argument per se, let me offer a clarifying of a few concepts: the first-one is the notion of a           {\em tree-presentation} and the second is the {\em finitization} of a tree through one of its tree-presentations. As we know, a formula-tree $\Phi$ is a potentially infinite object and we can represent it by an initial sub-tree having in some of its final nodes full sub-trees. Any such representation of $\Phi$ will be called a {\em tree-presentation} of $\Phi$; I consider it to be a finite object.

\smallskip It is also possible to ``finitize" entirely the presentation of a tree by replacing formulas in end nodes by relational symbols. For instance, consider the tree $\Phi(x,y)$ given by

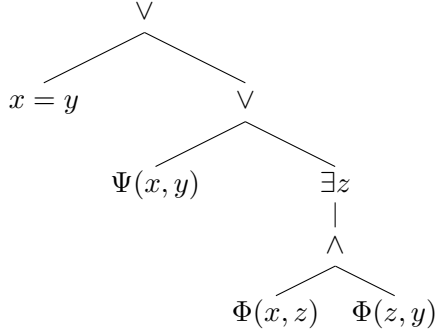
\Tree [.$\lor$ $x=y$ [.$\lor$ [.$\Psi(x,y)$ ] [.$\exists z$ [.$\land$ $\Phi(x,z)$ $\Phi(z,y)$ ] ] ] ]

\bigskip\noindent If we substitute $E^*$ for $\Phi$ and $E$ for $\Psi$ in the above formula, we obtain \cite[p.~217]{Halpern01}:
\begin{displaymath}
    E^*(x,y) \equiv x=y \lor E(x,y) \lor \exists z (E^*(x,z) \land E^*(z,y)).
\end{displaymath}

In the general case, since I am considering the language $\LNTplus$ as an infinite sequence of formula-trees, $\LNTplus=\langle\Phi_0,\ldots\Phi_n,\ldots\rangle$, what I do is to take an infinite sequence of relational symbols $\langle\textbf{R}_0,\ldots,\textbf{R}_n,\ldots\rangle$, each $\textbf{R}_n$ having the same arity of its corresponding $\Phi_n$, and substitute each $\textbf{R}_n$ in the end nodes of the presentation for $\Phi_n$.

\medskip
Now, I break the main proof into a few steps:
\begin{enumerate}
\item Let $\mathcal{TS}$ be an \textbf{NP}-complete problem. Then, by Fagin's theorem \ref{Fag} given above we may find a short existential second-order formula $\exists X\varphi(X)$ which may supposed to be concrete, such that $\mathcal{TS}=\{\mathcal{A}:\mathcal{A}\models\ulcorner\exists X\varphi\urcorner\}$.
\item $\mathcal{A}$ can be intuitively thought as constituting a finite database, with universe $|\mathcal{A}|$ of length $n$ encoded in number theory:
\begin{displaymath}
    \mathcal{TS}(\mathcal{A})\leftrightarrow\mathcal{A}\models\ulcorner\exists X\varphi\urcorner.
\end{displaymath}

\item Let $Tr_N^{(\mathcal{A})}$ be the truth definition for formulas of length $\leq lh(\exists X\varphi)=N$ in the database $\mathcal{A}$. It follows that
\begin{displaymath}
    \mathcal{TS}(\mathcal{A})\leftrightarrow Tr_N^{(\mathcal{A})}(\ulcorner\exists X\varphi\urcorner)
\end{displaymath}

\item Let's fix a database $\mathcal{\underline{A}}$ in $\mathcal{TS}$. It follows that  $\mathfrak{A}\models Tr_N^{(\mathcal{\underline{A}})}(\ulcorner\exists X\varphi\urcorner)$. From the definition of $Tr_N^{(\mathcal{\underline{A}})}$, we derive that $\mathfrak{A}\models\exists X Tr_N^{(\mathcal{\underline{A}})}(\ulcorner\varphi(X)\urcorner)[X/\ulcorner X\urcorner].$
    This comes from the fact that
    $$Tr_N^{(\mathcal{\underline{A}})}(\ulcorner\exists X\varphi\urcorner)\leftrightarrow(\exists X\varphi)^{(\mathcal{\underline{A}})} \leftrightarrow \exists X(X\subseteq |\mathcal{\underline{A}}|\land\varphi^{(\mathcal{\underline{A}})}).$$

\item
By the definition of truth in the Absolute $\mathfrak{A}$ for existential statements given at the end of the previous section, it follows that there exists a formula $\Phi(\mathcal{\underline{A}};v_0,\ldots v_{M-1})\in \LNTplus$ such that $\vdash_\mathbb{NT}^+ Tr_N^{(\mathcal{\underline{A}})}(\ulcorner\varphi(X)\urcorner)[\Phi/\ulcorner X\urcorner]$. It should be clear that if $S0^{n_0},\ldots,S0^{n_{M-1}}$ satisfy $\Phi(\mathcal{\underline{A}};v_0,\ldots,v_{M-1})$, i.e., $\vdash_\mathbb{NT}^+\Phi(\mathcal{\underline{A}};S0^{n_0},\ldots,S0^{n_{M-1}})$, then each $S0^{n_i}\in|\mathcal{\underline{A}}|$, for $0\leq i\leq{M-1}$.

\item Now, let's fix our attention to formula $\Phi(\mathcal{\underline{A}};v_0,\ldots v_{M-1})\in \LNTplus$. For this formula, there exists a pair of decisions $\delta^-,\delta^+$ such that for every finite sequence $\langle S0^{n_0},\ldots,S0^{n_{M-1}}\rangle\in |\mathcal{\underline{A}}|^M$, either $\delta^-(S0^{n_0},\ldots,S0^{n_{M-1}})$ is a proof of
$\lnot\Phi(\mathcal{\underline{A}};S0^{n_0},\ldots,S0^{n_{M-1}})$ or $\delta^+(S0^{n_0},\ldots,S0^{n_{M-1}})$ is a proof of $\Phi(\mathcal{\underline{A}};S0^{n_0},\ldots,S0^{n_{M-1}})$.

\item As we see, for each database $\mathcal{A}$ a relation $\textbf{R}_\mathcal{A}(v_0,\ldots,v_{M-1})$ is associated to the formula $\Phi$ by means of the pair of decisions $\delta^-,\delta^+$. For the fixed database $\mathcal{\underline{A}}$, each $M$-tuple $\langle S0^{n_0},\ldots,S0^{n_{M-1}}\rangle\in |\mathcal{\underline{A}}|^M$ is decided to satisfy $\textbf{R}_\mathcal{A}$ by a process of checking whether two sequences of formulas, $\delta^-(S0^{n_0},\ldots,S0^{n_{M-1}})$ and $\delta^+(S0^{n_0},\ldots,S0^{n_{M-1}})$, which one is a proof. This is a process accomplished in polynomial time.

\item Just for the sake of completeness, I'm going to state the Immerman-Vardi theorem \cite[p.~217]{Halpern01}:
\begin{Theorem}
({\sc Immerman-Vardi}). A problem is in polynomial time if and only if
it is describable in first-order logic with the addition of the least-fixed-point
operator. This is equivalent to being expressible by a first-order formula iterated
polynomially many times. Formally, $\mathrm{P = FO(LFP) = FO[n^{O(1)}]}$.
\end{Theorem}

\item We now use the stated theorem (of Immerman-Vardi \cite{Immerman86}, \cite{Vardi82}) to derive that the relation  $\textbf{R}_\mathcal{A}(v_0,\ldots,v_{M-1})$ is describable using first-order logic with the addition of the least-fixed-point operator, i.e., $\textbf{R}_\mathcal{A}(v_0,\ldots,v_{M-1})\in \text{FO(LFP)}.$

\item Finally, to check whether $\mathcal{\underline{A}}\in\mathcal{TS}$ one must verify whether $\mathcal{\underline{A}}\models\ulcorner\varphi\urcorner[\textbf{R}_\mathcal{\underline{A}}/\ulcorner X\urcorner]$, which is also a polynomial procedure, since $$\text{FO(FO(LFP))=FO(LFP)}.$$ \noindent Thus, $\mathcal{TS}$ is in \textbf{P}, and since $\mathcal{TS}$ is \textbf{NP}-complete, \textbf{P=NP}.$\square$

\end{enumerate}

\smallskip I would like to end this section with a few remarks:

\begin{itemize}
    \item In the eagerness to prove the previous result, I tried without success to use the Immerman-Vardi theorem \cite{Immerman86}, \cite{Vardi82}, by showing that the relation $\textbf{R}_\mathcal{A}(v_0,\ldots,v_{M-1})$ can be obtained by means of the fixed point operator. I think it to be an approach worth trying.

    \item I would like to make it plainly clear that I do not consider the argument above as a solution to the \PNP\  problem. This stems from the fact that Decision Theory has not been developed yet. But I haven't given up my hopes that computer scientists and logicians will in the future develop it, and fully confirm the above given argument. For the time being, I consider the \PNP\  problem as completely open.
\end{itemize}

    When asked about the \PNP\  problem, Professor Donald Knuth thus expressed himself \cite[Question 17]{Knuth14}:
        \begin{quotation}
            \textbf{Andrew Binstock, Dr. Dobb's:} ``At the ACM Turing Centennial in 2012, you stated that you were becoming convinced that \textbf{P=NP}. Would you be kind enough to explain your current thinking on this question, how you came to it, and whether this growing conviction came as a surprise to you?"

            \smallskip
            \noindent\textbf{Donald Knuth:} ``As you say, I've come to believe that \textbf{P=NP}, namely that there does exist an integer $M$ and an algorithm  that will solve every $n$-bit problem belonging to the class \textbf{NP} in $n^M$ elementary steps.

            Some of my reasoning is admittedly naïve: It's hard to believe that $\mathbf{P\not=NP}$ and that so many brilliant people have failed to discover why. On the other hand if you imagine a number $M$ that's finite but incredibly large --- like say the number $10\uparrow\uparrow\uparrow\uparrow3$ discussed in my paper on ``coping with finiteness" --- then there's a humongous number of possible algorithms that do $n^M$ bitwise or addition or shift operations on n given bits, and it's really hard to believe that all of those algorithms fail.

            My main point, however, is that I don't believe that the equality \textbf{P=NP} will turn out to be helpful even if it is proved, because such a proof will almost surely be nonconstructive. Although I think $M$ probably exists, I also think human beings will never know such a value. I even suspect that nobody will even know an upper bound on $M$.

            Mathematics is full of examples where something is proved to exist, yet the proof tells us nothing about how to find it. Knowledge of the mere existence of an algorithm is completely different from the knowledge of an actual algorithm."
        \end{quotation}

    I certainly agree with what Professor Knuth expresses, but I must mention that in the history of mathematics, there has occurred a strange, yet astonishing episode. The young David Hilbert solved the famous invariant problem, first by proving that a finite basis for the full system of invariants {\em does exist}, a strict existential proof\footnote{which lead Paul Gordan to declare that {\em this was Theology}.}. Afterwards, from the existence of a finite basis, Hilbert was able to produce a full constructive proof. Maybe an analogous fact will happen with  the \PNP\  problem.

In a paper published by the {\em Journal of Automated Reasoning} \cite{Thiele02}, Ruediger Thiele and Larry Wos discuss a presumable {\em twenty fourth} Hilbert problem, {\em hidden in his massive files in Germany.} The problem dealt with finding {\em simpler proofs and criteria for measuring simplicity.} Furthermore, the authors still discuss and analyse some methods of solution to the problem. Notwithstanding the fact that this  stems from Hilbert, this simplifying program for proofs might help finding a constructive proof from the existence of a polynomial algorithm for the travelling salesman problem.

\medskip
\section{Some Afterthoughts}
    To finish this article, I would like to make a few observations and raise a few questions which I think to be pertinent and related to second-order logic;  I would  rather shared them with you.

\medskip
\subsection{Second-order logic}
One of the most important problems connected to the ideas exposed in this article is the development of a deductive system for second-order logic, together with a completeness theorem, as it has been done for the first-order case. In order to accomplish it, we need to characterize the set-theoretical Absolute, which hopefully should be possible to do as it was done here, for number theory. The deductive system would need $\mathbb{ZF}$ amongst its axioms, and the completeness theorem should be carried out by taking all the ordinals, in correspondence with the first-order case which took only a set of newer constants. The second-order model for the given theory should be a full second-order model, which means that the set of theoretical predicates of the construing model would essentially coincide with the  power set of the model's universe. To do it, we need a full reflection principle on  the length and encode the theory in the language of second-order $\mathbb{ZF}$. A question that remains to be tackled with is whether the sentences derived in the deductive system are equivalent to existential second-order sentences, a fact that is valid in finite model theory under the hypothesis that \textbf{P=NP}. Finally, a second remark worth doing is that the formulas of the language $\mathcal{L}_\mathbb{ZF} ^+$ are sets, but the language itself is external to the set theoretical Absolute. This would arise difficulties in the development of the second-order deductive system, obliging us to use only finitely many or a recursively denumerable set of formulas out of $\mathcal{L}_\mathbb{ZF}^+$. This attempt could bring to light some of Silver glimpses into the nature of the set-theoretical Absolute. In the enlarged theory, I would not be surprised if in second-order $\mathbb{ZF}$ the power set axiom could be applied just once, to the set $\omega$.

\subsection{The continuum problem}
This problem appeared in the beginnings of set theory with Cantor. Essentially, what it is asked is the number of elements of the real line $\mathbb{R}$ measured by cardinals, or $\aleph$'s. Gödel showed that in the constructible world $L$, the power of the continuum is $\aleph_1$ (in $L$): ${2^{\aleph_0}}^L=\aleph_1^L$. He thus obtained his famous proof of {\em consistency} of the continuum hypothesis. Cohen, in his celebrated work, proved that the continuum problem is also independendent of the axioms of set theory. In his book, he thus expressed himself \cite[p.~151]{Cohen66}:

\begin{quotation}
``A point of view which the author feels may eventually come to be accepted is that CH is {\em obviously} false. The main reason one accepts the Axiom of Infinity is probably that we feel it absurd to think that the process of adding only one set at a time can exhaust the entire universe. Similarly with the higher axioms of infinity. Now $\aleph_1$ is the set of
countable ordinals and this is merely a special and the simplest way of
generating a higher cardinal, The set $C$ is, in contrast, generated by
a totally new and more powerful principle, namely the Power Set Axiom. It
is unreasonable to expect that any description of a larger cardinal which
attempts to build up that cardinal from the Replacement Axiom can ever reach $C$. Thus $C$ is greater than $\aleph_n$, $\aleph_\omega$, $\aleph_\alpha$, where $\alpha=\aleph_\omega$ etc.
This point of view regards
$C$ as an incredibly rich set
given to us by one bold new axiom, which can never be approached by any
piecemeal process of construction. Perhaps later generations will see
the problem more clearly and express themselves more eloquently."
\end{quotation}

Foreman, Magidor and Shelah \cite{Foreman88} as well as Woodin \cite{Woodin01, Woodin02}, by adjoining new principlles to $\mathbb{ZF}$, showed that the power of the  continuum is $\aleph_2$, a result much cherished by Gödel himself. In the late sixties, William Reinhardt mentioned a ``reflection principle on the width" and derived from it that the power of the continuum is very large, a weakly inaccessible cardinal. Nevertheless this result has never been published. The closest result  that I could find in the literature was Takeuti's Hypothesis 3 on the power set \cite[p.~440]{Takeuti68}, very akin to Reinhardt's principle. Let me suggest studying the possibility of defining the set theoretical Absolute together with a reflection principle on the width to derive that the power of the continuum is very large.

\subsection{Second-order computability and quantum computers}
The most outstanding theories of twentieth century physics are unquestionably Relativity theory and Quantum mechanics. Both theories have been married for almost a century, and they still do not communicate well. Relativity is a classical theory {\em par excellence} \cite{Andreka12}, and Quantum mechanics is what  it is \cite{Rovelli96}. Both theories are observer oriented.

\smallskip The great physicist David Deutsch modelled a universal Turing machine for quantum computers in the eighties \cite{Deutsch85}. He even stated a proposal that the Church-Turing thesis has a physical counterpart. Let me also suggest that classical physics is tied to first-order computability, linked to traditional Turing machines, descending from first-order proof theory. Whereas quantum Turing machines owed their origins to second-order computability, stemming from decision theory.

\smallskip After Cohen's discoveries, Professors Dana Scott and Bob Solovay developed independently, (and never published) the startling theory of Boolean valued models, and proved that it is in a certain sense ``equivalent"\ to Cohen's forcing. Since we relate in this section quantum theory to decision theory, and since quantum mechanics is probabilistic in a multi-varied number of interpretations, there might be a relationship in the development of decision theory with the Scott-Solovay  methods.

\subsection{The revival of Hilbert's formalist program}\label{SubsHFP}
Hilbert's formalist program for the foundations of mathematics suffered a major setback with Gödel's incompleteness theorems, despite the fact that Gödel himself did not agree \cite[p.~195]{Godel31} to it:
\begin{quotation}
    I wish to note expressly that Theorem XI (and the corresponding results for $M$ and $A$)\footnote{Theorem XI is Gödel's second incompleteness theorem, $M$ is the axiom system for set theory, $A$ is von Neumann's axiomatization of classical mathematics \cite{Neumann61}, and $P$ is {\em Principia Mathematicae}.} do not contradict Hilbert's formalist viewpoint. For this viewpoint presupposes only the existence of a consistency proof in which nothing but finitary means of proof is used, and it is conceivable that there exist finitary proofs that cannot be expressed in the formalism of $P$ (or of $M$ or $A$).
\end{quotation}

It is plausible that Hilbert's formalist program could be revived in our setting, since the formulas of $\LNTplus$ are potentially infinite but nevertheless finite objects, and therefore, Gödel's second incompleteness theorem no longer applies. If such a program is to be revived, we could start from and consider Finite Model Theory, even enriched with infinitary first-order formulas, as consistent.

\subsection{A global overview of classical mathematics}
By {\em classical mathematics}, at least in this article, I mean both versions of number theory, in first and second order logic. This is equivalent to saying that the hardcore of classical mathematics coincides with number theory and analysis. Nevertheless, such an assertion does not mean that set theory is being left out from classical mathematics, as we shall immediately see.

\smallskip
First of all, as mentioned in subsection \ref{SubsHFP}, and as a result of  bypassing Gödel's incompleteness theorems, the possibility of reviving Hilbert's formalistic program opens up again. I am firmly convinced that this will turn out to be the case. It is widely accepted that our physical universe precludes the actual infinite; therefore, actual infinite totalities may  be considered as ideal entities of the theories in which they appear. Allow me once again to quote Mrs.~Reid's book \cite[pp.~36--37.]{Reid96}:

\begin{quotation}
``When Klein went to Chicago for what was billed as an ``International Congress of Mathematicians'' to celebrate the founding of the University of Chicago, he took with him a paper by Hilbert in which that young man matter-of-factly summarized the history of invariant theory and his own part in it:

In the history of a mathematical theory the developmental stages are easily distinguished: the naive, the formal, and the critical. As for the theory of algebraic invariants, the first founders of it, Cayley and Sylvester, are together to be regarded as the representatives of the naive period: in the drawing up of the simplest invariant concepts and in the elegant applications to the solution of equations of the first degrees, they experienced the immediate joy of first discovery. The inventors and perfecters of the symbolic calculation, Clebsch and Gordan, are the champions of the second period. The critical period finds its expressions in the theorems I have listed above $\ldots$

The theorems he referred to were his own."
\end{quotation}

The foundations of mathematics go hand-in-hand with the foundations of set theory, and I see a mirroring between the development of the theory of invariants and that of set theory: the latter also had its naïve period, featuring Bolzano, Dedekind and Cantor; it then experienced its formal period, beginning with Frege, Russell, Zermelo, Hilbert, Gödel and Cohen, until the present; I conjecture that its critical period will be realized with the advent of second-order logic and the revival of Hilbert's program. Such a position could be underpinned through the remarkable results discovered by Harvey Friedman \cite{Friedman98}:

\begin{quotation}
``We present a coherent collection of finite mathematical theorems some of
which can only be proved by going well beyond the usual axioms for mathematics.
The proofs of these theorems illustrate in clear terms how one uses the
well studied higher infinities of abstract set theory called large cardinals in an
essential way in order to derive results in the context of the natural numbers.
The findings raise the specific issue of what constitutes a valid mathematical
proof and the general issue of objectivity in mathematics in a down to earth
way.

\noindent Large cardinal axioms, which go beyond the usual axioms for mathematics,
have been commonly used in abstract set theory since the 1960's (e.g., see
\cite{Scott61}, \cite{Martin89}). We believe that the results reported on here are the early
stages of an evolutionary process in which new axioms for mathematics will be
commonly used in an essential way in the more concrete parts of mathematics."
\end{quotation}

\noindent As we see, the down to earth proofs of the stated mathematical theorems require large cardinal axioms whose proof of consistency might perhaps be obtained in second-order number theory, as developed in this article.

\subsection{The Characteristica Universalis and Calculus Ratiocinator}
The Characteristica is considered as Leibniz's dream: a formal system to represent all knowledge together with a computational system to ``decide"\ all questions.
It could be epitomized in Leibniz's own words (translated by Bertrand Russell in \cite[p.~170]{Russell58}):

\begin{quotation}
``If controversies were to arise, there would
be no more need of disputation between two philosophers than
between two accountants. For it would suffice to take their pencils
in their hands, to sit down with their slates and say to each other
(with a friend as witness, if they liked): Let us calculate."
\end{quotation}
After Gödel's discovery of the incompleteness theorems, it became widely agreed that Leibniz's dream could never be materialized. This can also be seen from Carnap's notes taken in the Arkadencafé, in December 23, 1929, when Gödel presumably would have thus expressed himself \cite[p.~50]{Wang87}:

\begin{quotation}
``We admit as legitimate mathematics certain reflections on the grammar of a language that concerns the empirical. If one seeks to formalize such a mathematics, then with each
formalization there are problems, which one can understand and express in ordinary language, but cannot express in the given formalized language. It follows (Brouwer) that
mathematics is inexhaustible: one must always again draw afresh from the `fountain of intuition.' There is, therefore, no {\em Characteristica Universalis} for the {\em whole} of mathematics, and
no decision procedure for the whole mathematics. In each and every {\em closed language} there are only countably many expressions. The {\em continuum} appears only in `the whole of
mathematics.' $\ldots$ If we have {\em only one} language, and can only make `elucidations' about it, then these elucidations are inexhaustible, they always require some new intuition again,
tin contrast with these remarks, G's incompleteness theorem achieves the additional task of refining some of the ideas expressed here to get propositions {\em within} the given formal
system, which go beyond it relative to {\em provability}."
\end{quotation}

When Professor Abraham Robinson visited Brasília, in the early seventies, I asked him what were his thoughts about the {\em Characteristica Universalis}. His answer was brief: {\em ``This is a project for the future."} Well, the future is now: I believe that second order number theory, whether developed as in this article or not, together with the revival of Hilbert's formalistic program {\em is} the {\em Characteristica Universalis} with the {\em Calculus Ratiocinator}.

\subsection{Theology} Finally, I would like to remark that my readings of Gershom Scholem made me foresee a strong interaction between Theology and second-order logic, bridging the gap and promoting the unavoidable intertwining of Spirituality with the Foundations of Science.



\end{document}